\documentclass[a4paper,12pt]{amsart}

\usepackage{amsfonts} 
\usepackage{amssymb} 
\usepackage{amsmath,amsthm} 

\newtheorem{theoreme}{Theorem}[section] 
\newtheorem{lemma}[theoreme]{Lemma} 
 
\newtheorem{cor}[theoreme]{Corollary}

 \def \Rk {\ {\bf Remark.} } 
 
\def \sm {\setminus } 
\newcommand{\be}{\begin{enumerate}}  \newcommand{\ee}{\end{enumerate}} 
\newcommand{\bi}{\begin{itemize}}  \newcommand{\ei}{\end{itemize}} 
\newcommand{\bd}{\begin{description}}  \newcommand{\ed}{\end{description}}

\newcommand{\comment}[1]{}

\def \R {\mathbb{R}}

\numberwithin{equation}{section}  
\renewcommand{\phi}{\varphi} 
\renewcommand{\epsilon}{\varepsilon} 
 
\title{$L^p$-cohomology of negatively curved manifolds} 
\author{N. Yeganefar}
\thanks{This paper has been (partially) supported by the European Commission through the Research Training Network  HPRN-CT-1999-00118  "Geometric Analysis".} 
\date{January 23, 2004} 
\begin{document} 
\maketitle
\begin{abstract}
We compute the $L^p$-cohomology spaces of some negatively curved manifolds. We deal with two cases: manifolds with finite volume and sufficiently pinched negative curvature, and conformally compact manifolds.
\end{abstract} 
 \section{Introduction}
Let $(M,g)$ be a Riemannian manifold, and let $p\geq 1$ be a real number. We denote by $L^p(\Lambda ^*T^*M)$ (or $L^p$) the space of differential forms $\alpha$ with $|\alpha |^p$ integrable, and we consider the space $\Omega ^*_p(M)$ of elements of $L^p$ whose (weak) differential is also in $L^p$. The \emph{$L^p$-cohomology} of $(M,g)$, denoted by $H^*_p(M)$, is by definition the cohomology of the complex $(\Omega ^*_p(M),d)$:
 $$H^k_p(M)=\{ \alpha \in\Omega ^k_p(M)/d\alpha =0\}/d\Omega _p^{k-1}(M).$$ 
A difficulty which arises in the study of $L^p$-cohomology is the fact that the image $d\Omega _p^{k-1}(M)$ may not be closed in $L^p$, so that the topology of $H^k_p(M)$ may be complicated. One is then led to define \emph{reduced $L^p$-cohomology} by taking the quotient of $\{ \alpha \in\Omega ^k_p(M)/d\alpha =0\}$ by the closure of $d\Omega _p^{k-1}(M)$ in $L^p$. The reduced and unreduced $L^p$-cohomology spaces of a manifold are in general quite different, but if for example $H^k_p(M)$ is finite dimensional, then it is standard that they coincide. This happens for instance when $M$ is compact: it is known that in this case $H^k_p(M)$ is finite dimensional for all $k$, and in fact isomorphic to the $k^{th}$ de Rham cohomology group of $M$. 

For non compact $M$, we would like to know to what extent $H^k_p(M)$ reflects the topology or the geometry of the manifold. There exist already a lot of articles devoted to the study of this topic, see e.g. the works of Gol'dshtein, Kuzminov, Shvedov \cite{GKS1,GKS2}, Youssin \cite{Yo}, or Zucker \cite{Z2}, and the references therein. Of course, the case $p=2$ is of particular interest. For this case,  the proximity of  $L^2$-cohomology with the space of $L^2$ harmonic forms allows us to compute these (reduced or unreduced) $L^2$-cohomology spaces for large classes of manifolds, see the works of Zucker, Mazzeo and Mazzeo-Phillips \cite{Z1,M,MP} on some real or complex hyperbolic manifolds, of Dodziuk \cite{D} on rotationally symmetric manifolds, of Carron \cite{C} on flat manifolds, etc. 

In \cite{Y1,Y2}, we considered complete manifolds of finite volume and pinched negative curvature and identified the $L^2$-cohomology spaces with topologically defined groups (with a sharp pinching condition for real manifolds, and without any pinching assumption for K\"ahler manifolds). Here, our main goal is to deal more generally with the $L^p$-cohomology of these manifolds. Our main result is
\begin{theoreme}\label{coho}
Let $(M,g)$ be a complete $n-$dimensional manifold of finite volume and pinched negative sectional curvature $K$, i.e. there exists a constant $a>0$ such that $-1\leq K\leq -a^2<0$. Assume that $p\geq 1$ is a real number, and $k$ an integer such that $(k-2)(p-1)a-(n-k+1)>0$. Then we have the isomorphism $H^k_p(M)\simeq H^k_c(M)$, where $H^k_c(M)$ denotes the compactly supported cohomology of $M$.
\end{theoreme}
Thus, if the curvature is not too far from $-1$, we can compute the $L^p$-cohomology spaces in terms of the topology of the manifold. Roughly speaking, to prove this theorem we will proceed as follows. First, there is an exact sequence (see \cite{GKS1}) which relates the compactly supported cohomology of a bounded subset $D\subset M$ to the $L^p$-cohomology of $M$ and $M\sm D$. We will then use an idea developed by Pansu \cite{P} in order to prove that the $L^p$-cohomology of $M\sm D$ vanishes, and finally conclude by invoking the exact sequence.

 We will also see that our method works for other cases as well. For example, it can be applied to conformally compact manifolds (the definition will be recalled in the last section). The $L^2$-cohomology of these manifolds was computed by R. Mazzeo \cite{M} (see also \cite{Y2} for a simpler proof). We then state our next result, which is also a consequence of \cite{GKS1} (see the remark at the end of section 4):
\begin{theoreme}\label{cc}
Let $(M^n,g)$ be a conformally compact $n$-dimensional manifold. Assume that $p\geq 1$ is a real number and $k$ an integer such that $n-1-kp>0$. Then we have the isomorphism $H^k_p(M)\simeq H^k_c(M)$. 
\end{theoreme}
The paper is organized as follows: first, in the framework of negatively curved manifolds of finite volume, we define a homotopy operator acting on differential forms and study its main properties. Then we use this operator to prove Theorem \ref{coho} and deduce from this proof two interesting corollaries. We finally prove Theorem \ref{cc}.\\

{\it Acknowledgements.} I would like to thank my PhD supervisor G. Carron for many helpful discussions, hints and comments.  I  warmly thank J. Br\"uning and his team, especially G. Marinescu, for their hospitality at the Humboldt Universit\"at (Berlin), while part of this work was written.
\section{The homotopy operator}
In this section, we will introduce a "homotopy operator" acting on differential forms, which will be used to prove Theorem \ref{coho}. We are inspired by Pansu's work \cite{P}. From now on  $(M,g)$ will be  a complete $n-$dimensional manifold of finite volume and pinched negative sectional curvature $K$, i.e. there exists a constant $a>0$ such that 
$$-1\leq K\leq -a^2<0.$$
\subsection{Geometry of such manifolds}
Here, we briefly recall some standard facts about the topology and geometry of these manifolds (see
\cite{E} and \cite{HI}). First, $M$ has a finite number of ends, and one has $M=M_0\cup E_i$, where $M_0$ is a compact manifold with boundary, and the $\partial E_i$'s are the components of $\partial M_0$. To each ray of $E_i$, we can associate a Busemann function $r_i$ which is
 a priori only $C^2$-smooth. Two such functions are equal up to an additive
 constant. $E_i$ is $C^2-$diffeomorphic to $\R^+\times \partial E_i$. Moreover, the slices $\{ t\}\times \partial E_i$ are the level sets of a Busemann function. Finally, the metric on each end $E_i$ has the following form
$$g=dr_i^2+h_{r_i},$$
where $h_{r_i}$ is a family of metrics on the compact manifold $\partial E_i$, and satisfies  $e^{-r_i}h_0\leq h_{r_i}\leq e^{-ar_i}h_0.$

\subsection{Two technical lemmas}
In order to prove Theorem \ref{coho}, we  need to prove first two results which are the analogs of \cite[Propositions 8 and 10]{P}. To simplify notation, we assume that our manifold $M$ has only one end, and we choose an associated Busemann function $r$ such that outside the bounded open subset $D=\{ r<0 \}$, we have $M\sm D\simeq [0,\infty )\times \partial D$. Let $\nabla r$ denote the gradient of $r$, with flow $\phi _t$. On $M\sm D$ the flow $\phi _t$ is just  $t$-translation on the first factor $[0,\infty )$; more precisely, if $x=(r_0,\theta _0)$ is a point in $M\sm D$, then we simply have $$\phi _t(r_0,\theta _0)=(r_0+t,\theta _0).$$
Following \cite{P}, we now introduce the operator $B$ which acts on a differential form $\alpha$ by
\begin{equation}\label{operateurB}
B\alpha = \int_0^\infty \phi _t^*(i_{\nabla r}\alpha)\,dt,
\end{equation}
where $i_{\nabla r}\alpha$ denotes interior product. In the next two lemmas, we study the main properties of this operator. First, we have
\begin{lemma}\label{contraction}
Denote by $Jac(\phi_t)$ the Jacobian of $\phi _t$. For an integer $j$, set $$\eta =(j-1)(p-1)a-(n-j).$$ Then for every $j-$forme $\alpha$ defined on $M\sm D$, and for all $x$ in $M\sm D$, we have 
$$ |\phi _t^*(i_{\nabla r}\alpha)(x)|^p\leq e^{-\eta t}|\alpha (\phi _t(x))|^p Jac(\phi _t)(x).$$
In particular, if $\eta >0$, the operator $B$ is well defined and bounded on $L^p(\Lambda ^jT^*(M\sm D))$.
\end{lemma}
We don't give the proof. The one in \cite[Proposition 8]{P} can be applied word by word to our case: roughly speaking, we have to estimate the derivative of $p\log (|\phi _t^*(i_{\nabla r}\alpha)(x)|/|\alpha (\phi _t(x))|)- \log (Jac(\phi _t)(x))$ with respect to $t$, and to do this we have to estimate the principal curvatures of the level sets $\{ t\}\times \partial D$  of the Busemann function. The only difference is that these principal curvatures are between $-1$ et $-a$, and not between $a$ et $1$ as in the case of negatively curved simply connected manifolds considered in \cite[Proposition 8]{P}.

Next, we have
\begin{lemma}\label{homotopy}
Let $\alpha$ be a $j-$form in $\Omega ^j_p(M\sm D)$. Suppose that $(j-1)(p-1)a-(n-j)>0$. Then we have the following homotopy formula
 $$dB\alpha +Bd\alpha =-\alpha .$$
\end{lemma}
\begin{proof}
We first assume that $\alpha$ is smooth and with compact support in $M\sm D$ (this support can meet the boundary $\partial D$). By Cartan formula, we get  
\begin{eqnarray*}
d\phi _t^*i_{\nabla r}\alpha +\phi _t^*i_{\nabla r}d\alpha &=& \phi _t^*(di_{\nabla r}\alpha +i_{\nabla r}d\alpha )\\           &=&\phi _t^* L_{\nabla r}\alpha \\
              &=&\frac{\partial}{\partial t}(\phi _t^*\alpha ),
\end{eqnarray*}
where $L_{\nabla r}$ is Lie derivative with respect to the vector field $\nabla r$. It follows that
\begin{eqnarray*}
dB\alpha +Bd\alpha &=&\int _0^\infty \frac{\partial}{\partial t}(\phi _t^*\alpha )\,dt \\
                   &=&\lim_{t\to\infty}\phi _t^*\alpha-\alpha .
\end{eqnarray*}
By the properties of the flow $\phi _t$, we have $\lim_{t\to\infty}\phi _t^*\alpha =0$ because $\alpha$ has compact support. Therefore the homotopy formula holds for smooth compactly supported forms. Now, the assumption on the numbers $a$, $p$ and $j$ implies that the operator $B$ is bounded on $j$ and $(j+1)$-forms which are $L^p$ on $M\sm D$ (see Lemma \ref{contraction}).  But the space of smooth $j-$forms with compact support on $M\sm D$ is dense in $\Omega ^j_p(M\sm D)$, so that the formula is valid on $\Omega ^j_p(M\sm D)$.
\end{proof}

\section{$L^p$-cohomology of negatively curved manifolds}
\subsection{Proof of Theorem \ref{coho}}
We keep our previous notations. By \cite[Lemma 10]{GKS1}, we have the following exact sequence:
$$H^{k-1}_p(M\sm D)\buildrel b\over \to H^k_c(D)\buildrel e\over \to H^k_p(M)\buildrel r\over \to H^k_p(M\sm D),$$
where $b$ is the coboundary operator, $e$ is extension by zero, and $r$ restriction.
Thus, to prove Theorem \ref{coho}, it is enough to show that the $L^p$-cohomology spaces at infinity $H^j_p(M\sm D)$ vanish for $j=k-1,k$.

Let $[\alpha ]$ be an element in $H^j_p(M\sm D)$, with $j=k-1$ or $k$, and let $\alpha \in \Omega ^j_p(M\sm D)$ be a representative of this class. The hypothesis $(k-2)(p-1)a-(n-k+1)>0$ implies that the assumptions of Lemma \ref{homotopy} are satisfied for $j=k-1,k$. As $\alpha$ is closed, we therefore get 
$$\alpha =d(-B\alpha).$$
But $\alpha$ is in $\Omega ^j_p(M\sm D)$, so that $B\alpha$ is also in $\Omega ^j_p(M\sm D)$ by Lemma \ref{contraction}; hence $\alpha$ is zero in $H^j_p(M\sm D)$.

\subsection{Further results}
If we only assume that $(k-1)(p-1)a-(n-k)>0$, then the proof of Theorem \ref{coho} still gives the vanishing of the $L^p$-cohomology group at infinity  $H^k_p(M\sm D)$. This shows, via the exact sequence, that there is a surjection of $H^k_c(D)$ onto  $H^k_p(M)$, so that $H^k_p(M)$ is in particular finite dimensional. This has the following well-known consequence:
\begin{cor}
We keep the same notations as in Theorem \ref{coho}, and assume that $(k-1)(p-1)a-(n-k)>0$. Then the range of $d :\Omega ^{k-1}_p(M)\to L^p(\Lambda ^kT^*M)$ is closed.
\end{cor}
\Rk Under this hypothesis $(k-1)(p-1)a-(n-k)>0$, which is weaker than the one in Theorem \ref{coho}, there is probably still an isomorphism $H^k_p(M)\simeq H^k_c(M)$. This is true if $p=2$ (see \cite[Proposition 5.1]{Y2}).

By a duality argument, we can get another corollary from Theorem \ref{coho}. Namely, for any real number $p>1$, denote by $q$ the conjugate exponent, i.e. $1/p+1/q=1$. If we assume that $M$ is orientable, then the bilinear form $H^k_p\times H^{n-k}_q \to\R$ given by 
$$([\alpha ],[\beta ])\mapsto \int _M\alpha \wedge \beta $$ is well defined and is non degenerate (see e.g. \cite[Lemme 81]{P}). We therefore obtain
\begin{cor}
Let $M$ be a complete orientable $n-$dimensional Riemannian manifold with finite volume and pinched negative curvature $-1\leq K\leq -a^2<0$. If $p>1$ and $k$ satisfy $(k-2)(p-1)a-(n-k+1)>0$, then we have the isomorphism $H^{n-k}_q(M)\simeq H^{n-k}(M)$.
\end{cor}

\section{Conformally compact manifolds}
The method of proof of Theorem \ref{coho} can be extended to other situations, too. As an illustration, we will consider the case of conformally compact manifolds. Let us recall some basic facts about the geometry of such manifolds.

Let $\overline{M}$ be a compact manifold with boundary, equipped with a Riemannian metric $\overline{g}$ which is smooth up to the boundary. Let $y : \overline{M}\to \R ^{+}$ be a nonnegative smooth defining function for the boundary $\partial M$: $\partial M=y^{-1}(0)$ and $dy\neq 0$ along $\partial M$. We endow $M$ with the metric $g=\overline{g}/y^{2}$; this metric is complete, and we say that $(M,g)$ is a \emph{conformally compact} manifold. The typical example of such a situation is the ball model of the real hyperbolic space. With these preliminaries in mind, we can now prove Theorem \ref{cc}.\\
{\it \bf Proof of Theorem \ref{cc}.}\\
It is well-known (see e.g. the proof of \cite[Corollary 6.2]{Y2}) that outside a compact subset $D$, $M$ is quasi-isometric to a warped product $([0,\infty )\times \partial M, dr^2+e^{2r}h)$, where $h$ is a metric on $\partial M$ which does not depend on $r$. As $L^p$-cohomology is by definition invariant under quasi-isometries, it will be enough to consider this warped product case. We follow the same line of reasoning as in the proof of Theorem \ref{coho}. Thus, we define the operator $B$ by the formula \ref{operateurB}. The principal curvatures of the level sets $\{ t\} \times \partial M$ of the function $r$ are all equal to $1$. Hence the proof of \cite[Proposition 8]{P} shows that for every $j-$form $\alpha$ defined on $M\sm D$, we have 
$$|\phi _t^*(i_{\nabla r}\alpha)(x)|^p\leq e^{-(n-1-(j-1)p)t}|\alpha (\phi _t(x))|^p Jac(\phi _t)(x).$$
Therefore, if $n-1-(j-1)p>0$, then $B$ is bounded on $L^p(\Lambda ^jT^*(M\sm D))$. Now, if $n-1-kp>0$, $B$ is bounded on $(k-1)$ and $k-$forms, so that we have the homotopy formula of Lemma \ref{homotopy} for elements of $\Omega ^{k-1}_p(M\sm D)$  and $\Omega ^{k}_p(M\sm D)$. From this, we deduce the vanishing of the $L^p-$cohomology spaces at infinity in degrees $k-1$ and $k$, and finally complete the proof by using the exact sequence of \cite{GKS1}.

\Rk As pointed out in the introduction, Theorem \ref{cc} is also a consequence of \cite{GKS1}. Namely, for warped products, Gol'dshtein, Kuz'minov and Shvedov get vanishing results for some  $L^p$-cohomology spaces at infinity , and we can then use the exact sequence.


{\small \textsc{Universit\'e de Nantes, D\'epartement de Math\'ematiques, 2 rue de la Houssini\`ere, BP 92208, 44322 Nantes cedex 03, France\\}}
e-mail: nader.yeganefar@math.univ-nantes.fr
\end{document}